\documentclass[11pt]{article}
\usepackage{mathrsfs}
\usepackage{amsfonts} \makeatother
\usepackage{amsfonts,amsmath,amssymb,mathrsfs}

\numberwithin{equation}{section}
\newtheorem{theorem}{Theorem}[section]
\newtheorem{lemma}[theorem]{Lemma}

\newtheorem{remark}[theorem]{Remark}

\newcommand{\be}{\begin{equation}}
\newcommand{\ee}{\end{equation}}
\newcommand\bes{\begin{eqnarray}} \newcommand\ees{\end{eqnarray}}
\newcommand{\bess}{\begin{eqnarray*}}
\newcommand{\eess}{\end{eqnarray*}}

\usepackage{amssymb}
\begin{document}

\begin{center} {\bf\Large Multiple positive solutions for Kirchhoff type problems involving}
\\[2mm]
 {\bf\Large  concave and critical nonlinearities in $\mathbb{R}^3$}
\\[4mm]
{\large \ \ Xiaofei Cao$^{\dag}$,
\ \ Junxiang Xu$^{\dag,}$\footnote{Corresponding author.

~~\textit{E-mail\ addresses}: caoxiaofei258@126.com(X. Cao),\
xujun@seu.edu.cn (J. Xu),\ wangmath2011@126.com (J. Wang).\ }} ,\ \
\ \ Jun Wang$^{\ddag}$  \\[2mm]
{$^{\dag}$ Department of Mathematics,  Southeast University,

 Nanjing
210096, P.R. China.

$^{\ddag}$ Faculty of Science, Jiangsu University, Zhenjiang,

Jiangsu, 212013, P.R. China.}\\[2mm]
\end{center}

\setlength{\baselineskip}{16pt}

\begin{quote}
\noindent {\bf Abstract:} In this paper, we consider the multiplicity of solutions for a class of Kirchhoff type problems with sub-linear and critical terms on an unbounded domain. With the aid of Ekeland's variational principle and the concentration compactness principle we prove that the Kirchhoff problem has at least two solutions.

\noindent {\bf Keywords}: Kirchhoff type problem, the concentration compactness principle, variational method.

\noindent{\bf {MSC(2010)}}: 35A01, 35A15.
\end{quote}

\section{Introduction and main results}

\setcounter{section}{1} \setcounter{equation}{0}

\noindent
 This paper concerns the multiplicity of solutions for the following Kirchhoff type problem
\begin{equation}\label{a1}
\begin{cases}
-(a+b\int_{\mathbb{R}^{3}}|\nabla u|^{2}dx)\triangle u=\lambda f(x)|u|^{q-2}u+u^{5},~~~~x\in\mathbb{R}^{3},\\
u\in D^{1,2}(\mathbb{R}^{3}),\\
\end{cases}
\end{equation}
where $a, b$ are positive constants, $f(x)$ is a continuous function, $1\leq q\leq2$ or $1\leq q<2$ .

It is well known that Kirchhoff type problems are proposed by Kirchhoff in 1883 \cite{GGKKM} as an extension of the classical D'Alembert's wave equation
for free vibration of elastic strings. Such problems are often viewed as nonlocal because the presence of the integral term $\int|\nabla u|^{2}dx$. This phenomenon causes some mathematical difficulties making the study of such problems particularly interesting. The case of Kirchhoff problems where the nonlinear term is super-linear has been investigated in the last decades by many authors, for example \cite{BCXW,DGW,GMFG,XHWZ,XHWWZ,JJXW,WLXH,YLFL,JSSL,JWLX,ZTZKP} and references therein. Here, we are interested in the case of Kirchhoff problems where the nonlinear term is sub-linear and critical.

For nonlinear elliptic problems, Chabrowski and Drabek \cite{Chabrowski} considered the following nonlinear elliptic problem:
\begin{equation}\label{a2}
-\triangle u+V(x)u=\epsilon h(x)u^{q}+u^{2^{\ast}-1} \quad \text{in}\quad\mathbb{R}^N,
\end{equation}
where $\epsilon>0$ is a parameter, $1<q<2$ and $2^{\ast}=\frac{2N}{N-2}, N\geq3$, is the critical Sobolev exponent. Under the assumptions of $h$ is a nonnegative and nonzero function in $L^{r}(\mathbb{R}^N)\cap C(\mathbb{R}^N)$, where $r=\frac{2^{\ast}}{2^{\ast}-q-1}$, they obtained that \eqref{a2} has at least two nonnegative solutions by applying various variational principle.

For Kirchhoff problems, Fan \cite{Fan} investigated the existence of multiple positive solutions to the following Kirchhoff type problem:
\begin{equation}\label{a3}
\begin{cases}
-(a+b\int_{\Omega}|\nabla u|^{2}dx)\triangle u=d(x)u^{k-1}u+g(x)u^{5},\ &x\in\Omega,\\
u=0,\ &x\in\partial\Omega,\\
\end{cases}
\end{equation}
where $a, b>0, 4<k<6$, $\Omega$ is a smooth bounded domain in $\mathbb{R}^3$ and $d(x), g(x)$ are positive and continuous functions. By introducing suitable  conditions on $d(x), g(x)$, they proved that there exists $\Lambda_{\delta}$ such that if $|f|_{L^{q^{\ast}}}<\Lambda_{\delta}$, \eqref{a3} has at least $cat_{M_{\delta}}(M)$ distinct positive solutions, where $q^{\ast}=\frac{6}{6-q}$ and $cat$ mens the Ljusternik-Schnirelmann category(see \cite{MWM}).

Liu et. al.\cite{Liu} considered the following nonlinear Kirchhoff type equation
\begin{equation}\label{a6}
\begin{cases}
-(a+b\int_{\mathbb{R}^{N}}|\nabla u|^{2}dx)\triangle u=u^{2^{\ast}-2}u+\mu h(x),~~~~x\in\mathbb{R}^{N},\\
u\in D^{1,2}(\mathbb{R}^{N}),\\
\end{cases}
\end{equation}
where $a\geq0, b>0, N\geq3, 2^{\ast}=\frac{2N}{N-2}, \mu\geq0$ and $h\in L^{\frac{2^{\ast}}{2^{\ast}-1}}(\mathbb{R}^{N}\backslash \{0\})$ is nonnegative. Under some assumptions on $a, b$ and $\mu$, they obtained the existence of two positive solutions for Eq. \eqref{a6}.

Sun \cite{JSSL} studied the following Kirchhoff type problem with critical exponent
\begin{equation}\label{a5}
\begin{cases}
-(a+b\int_{\Omega}|\nabla u|^{2}dx)\triangle u=\lambda |u|^{q-2}u+u^{5},\ &x\in\Omega,\\
u=0,\ &x\in\partial\Omega,\\
\end{cases}
\end{equation}
where $a, b, \lambda>0, 1<q<2$, $\Omega$ is a smooth bounded domain in $\mathbb{R}^3$. They showed that there exists a positive constant $T(a)$ depending on $a$ such that for each $a>0$ and $0<\lambda<T(a)$, \eqref{a5} has at least one positive solution.

Xie et. al. \cite{Xie} considered the following Kirchhoff type problems
\begin{equation}\label{a4}
\begin{cases}
-(a+b\int_{\mathbb{R}^{3}}|\nabla u|^{2}dx)\triangle u+V(x)u=u^{5},~~~~x\in\mathbb{R}^{3},\\
u\in D^{1,2}(\mathbb{R}^{3}),\\
\end{cases}
\end{equation}
where $a,b>0$ and $V\in L^{\frac{3}{2}}(\mathbb{R}^{3})$ is a given nonnegative function. If $|V|_{\frac{3}{2}}$ is suitable small, they proved that \eqref{a4} has at least one bound state solution.

Motivated by above papers, we consider the Kirchhoff problem \eqref{a1} with concave and critical nonlinearities on the whole space $\mathbb{R}^{3}$. 
To the best of our knowledge, there are few papers which deal with this type of Kirchhoff problem \eqref{a1}. The main difficulty is how
to estimate the energy and recover the compactness
because the nonlinearity is the combination of the concave and critical terms. By the method of Mountain Pass Theorem and the concentration compactness principle, we obtain \eqref{a1} has at least two different solutions with their energies having different signs.

\begin{theorem}\label{thm 1.1}
Assume that in the problem \eqref{a1}, $a, b$ are positive constants, $1\leq q\leq2$ and $f(x)$ is a function that can change sign and with property
 \begin{enumerate}
  \item [$(f)$]$f(x)\in C(\mathbb{R}^{3})\cap L^{q^{\ast}}(\mathbb{R}^{3})$, where $q^{\ast}=\frac{6}{6-q}$.
\end{enumerate}
Then there exists $\lambda_{1}>0$ such that if $\lambda\in(0,\lambda_{1})$, the problem \eqref{a1} has one positive solutions which has a negative energy.
\end{theorem}

\begin{theorem}\label{thm 1.2}
Assume that in the problem \eqref{a1}, $a, b$ are positive constants, $1\leq q<2$ and $f(x)$ is a nonnegative function with property $f$.
Then there exists $0<\lambda_{2}\leq\lambda_{1}$ such that if $\lambda\in(0,\lambda_{2})$, the problem \eqref{a1} has two positive solutions, one of which has a positive energy and the other a negative energy.
\end{theorem}

\begin{remark}\label{rema 1.3}
Compare with Fan \cite{Fan}, we consider the Kirchhoff problem with the nonlinear term is sub-linear(or linear) and critical in the while space $\mathbb{R}^{3}$ i.e. $1\leq q<2$ in Eq.  \eqref{a1}, while he investigated the Kirchhoff problem with the nonlinear term is super-triple and critical in a smooth bounded domain in $\mathbb{R}^{3}$ i.e. $4<k<6$ in Eq. \eqref{a3}. Compare with Liu et. al.\cite{Liu}, if $q=1$, $f(x)$ is a nonnegative function and $f(x)\in C(\mathbb{R}^{3})\cap L^{\frac{2^{\ast}}{2^{\ast}-1}}(\mathbb{R}^{N}\backslash \{0\})$ in Eq. \eqref{a1}, then
\begin{equation*}
\begin{cases}
-(a+b\int_{\mathbb{R}^{3}}|\nabla u|^{2}dx)\triangle u=\lambda f(x)+u^{5},~~~~x\in\mathbb{R}^{3},\\
u\in D^{1,2}(\mathbb{R}^{3}),\\
\end{cases}
\end{equation*}
which is the same as Eq. \eqref{a6} when $N=3$. Compare with Sun \cite{JSSL}, from Theorem \ref{thm 1.2}, we obtain the existence of two positive solutions for Eq. \eqref{a1} in the whole space $\mathbb{R}^{3}$,  while he obtained the existence of one positive solution for Eq. \eqref{a5} in a smooth bounded domain in $\mathbb{R}^{3}$. Compare with Xie et. al. \cite{Xie}, if $\lambda=1$, $q=2$, $f(x)$ is a negative function and $f(x)\in C(\mathbb{R}^{3})\cap L^{\frac{3}{2}}(\mathbb{R}^{3})$ in Eq. \eqref{a1}, then
\begin{equation*}
\begin{cases}
-(a+b\int_{\mathbb{R}^{3}}|\nabla u|^{2}dx)\triangle u+f(x)u=u^{5},~~~~x\in\mathbb{R}^{3},\\
u\in D^{1,2}(\mathbb{R}^{3}),\\
\end{cases}
\end{equation*}
which is the same as \eqref{a4}. From Theorem 1.1, we can see that there exists $\sigma>0$ such that if $|f|_{\frac{3}{2}}\in(0, \sigma)$, the problem \eqref{a1} has one positive solutions which has a negative energy.
\end{remark}

Throughout this paper, we make use of the following notations:

$\bullet$ $\rightarrow$(respectively, $\rightharpoonup$) denotes strong (respectively, weak) convergence;

$\bullet$ $|\cdot|_{p}$ denotes the norm of $L^{p}(\mathbb{R}^{N})$;

$\bullet$ $D^{1,2}(\mathbb{R}^{3})$ denotes the usual Sobolev space equipped with the norm $\|u\|= (\int_{\mathbb R^{3}}|\triangledown u|^{2}dx)^{\frac{1}{2}}$.

$\bullet$ $X^{\ast}$ denotes the dual space of $X$;

$\bullet$ $B_{\alpha}:=\{u\in D^{1,2}(\mathbb{R}^{3}): \|u\|=\alpha\}$ and $\overline{B}_{\alpha}:=\{u\in D^{1,2}(\mathbb{R}^{3}): \|u\|\leq\alpha\}$.

$\bullet$ $C, C_{1}, C_{2},\ldots$ denote various positive constants, which may vary from line to line.

This  paper is  organized  as follows: Section 2 is dedicated to the
abstract framework and some preliminary results. Sections 3 and 4 are dedicate to the proofs of Theorems \ref{thm 1.1} and \ref{thm 1.2}, respectively.

\section{Preliminaries}

The energy functional corresponding to \eqref{a1} is defined on $D^{1,2}(\mathbb R^{3})$ by
\begin{equation}\label{b1}
I(u)=\frac{a}{2}\| u\|^{2}+\frac{b}{4}\| u\|^{4}-\frac{1}{q}\lambda\int_{\mathbb R^{3}}f(x)|u|^{q}dx-\frac{1}{6}\int_{\mathbb R^{3}}|u|^{6}dx.
\end{equation}
It is well known that a weak solution of problem \eqref{a1} is a critical point of the functional $I$. In the following, we are devoted to finding critical points of $I$.

\begin{lemma}\label{lem 2.1}(\cite{MWM})
Suppose the hypothesis $(f)$ holds, the function $\varphi: D^{1,2}(\mathbb R^{3})\mapsto\mathbb R$ defined by
\begin{equation*}
\varphi(u)=\int_{\mathbb R^{3}}f(x)|u|^{q}dx
\end{equation*}
is weakly continuous. Moreover, $\varphi$ is continuously differentiable with derivative $\varphi': D^{1,2}(\mathbb R^{3})\mapsto (D^{1,2}(\mathbb R^{3}))^{\ast}$ given by
\begin{equation*}
\langle \varphi^{'}(u),v\rangle=q\int_{\mathbb R^{3}}f(x)|u|^{q-2}uvdx.
\end{equation*}
\end{lemma}

Obviously, the functional $I\in C^{1}(D^{1,2}(\mathbb R^{3}),\mathbb{R})$ and for any $u,v\in D^{1,2}(\mathbb R^{3})$,
\begin{equation*}
\langle I^{'}(u),v\rangle=(a+b\| u\|^{2})\int_{\mathbb R^{3}}\triangledown u\triangledown vdx-\lambda\int_{\mathbb R^{3}}f(x)|u|^{q-2}uvdx-\int_{\mathbb R^{3}}|u|^{4}uvdx.
\end{equation*}

\begin{lemma}\label{lem 2.2}(\cite{MWM}) The embedding $D^{1,2}(\mathbb R^{3})\hookrightarrow L^{6}(\mathbb R^{3})$ is continuous. Denote by $S$ the best Sobolev constant, which is given by
\begin{equation}\label{b2}
S:=\inf_{u\in D^{1,2}(\mathbb R^{3})\backslash\{0\}}\frac{\int_{\mathbb R^{3}}|\triangledown u|^{2}dx}{(\int_{\mathbb R^{3}}|u|^{6}dx)^{\frac{1}{3}}}.
\end{equation}
Moreover, $S$ is achieved by the function
\begin{equation}\label{b3}
U_{\epsilon}(x)=\frac{(3\epsilon^{2})\frac{1}{4}}{(\epsilon^{2}+|x|^{2})^{\frac{1}{2}}}.
\end{equation}
and
\begin{equation}\label{b4}
\int_{\mathbb R^{3}}|\triangledown U_{\epsilon}(x)|^{2}dx=\int_{\mathbb R^{3}}|U_{\epsilon}(x)|^{6}dx=S^{\frac{3}{2}}.
\end{equation}
\end{lemma}

\section {Proof of Theorem 1.1}
In this section, we are devoted to the proof of Theorem \ref{thm 1.1}, so we suppose that the assumptions of Theorem \ref{thm 1.1} hold throughout this section.
\begin{lemma}\label{lem 3.1}
\begin{enumerate}
\item [$(i)$] There exists $\lambda_{1}>0$ such that if $\lambda\in(0,\lambda_{1})$, then there exist $\alpha>0$ and $\rho>0$ such that
\begin{equation*}
I(u)|_{B_{\alpha}}\geq\rho>0.
\end{equation*}
\item [$(ii)$] There is $u_{0}\in \overline{B}_{\alpha}$ such that $I(u_{0})<0$.
\end{enumerate}
\end{lemma}

\textbf{proof}(i) By the assumption $(f)$, the H$\ddot{o}$lder inequality and Lemma \ref{lem 2.1}, we have
\begin{equation}\label{c1}
\begin{split}
I(u)&\geq\frac{a}{2}\|u\|^{2}+\frac{b}{4}\|u\|^{4}
-\frac{\lambda}{q}|f|_{q^{*}}|u|_{6}^{q}-\frac{1}{6}|u|^{6}_{6}\\
&\geq\frac{a}{2}\|u\|^{2}+\frac{b}{4}\|u\|^{4}-\frac{\lambda}{qS^{\frac{q}{2}}}|f|_{q^{*}}\|u\|^{q}-\frac{1}{6S^{3}}\|u\|^{6}\\
&=\|u\|^{q}\big(\frac{b}{4}\|u\|^{4-q}-\frac{\lambda}{qS^{\frac{q}{2}}}|f|_{q^{*}}-\frac{1}{6S^{3}}\|u\|^{6-q}\big).
\end{split}
\end{equation}

Set $l(t)=\frac{b}{4}t^{4-q}-\frac{1}{6S^{3}}t^{6-q}$ for $t>0.$ Direct calculations yield that
\begin{equation*}
\max_{t>0}l(t)=l(\alpha)=(\frac{3bS^{3}(4-q)}{2(6-q)})^{\frac{4-q}{2}}\cdot\frac{b}{2(6-q)}:=C_{p,q},
\end{equation*}
where $\alpha=(\frac{3bS^{3}(4-q)}{2(6-q)})^{\frac{1}{2}}$. Then it follows from \eqref{c1} that, if $\lambda<\lambda_{1}$, $I(u)|_{B_{\alpha}}\geq\rho>0$, where $\lambda_{1}=qS_{2}^{\frac{q}{2}}C_{p,q}\cdot\frac{1}{|f|_{q^{*}}}$ and $\rho=\alpha^{q}(l(\alpha)-\frac{\lambda}{q}|f|_{q^{\ast}}S_{2}^{-\frac{q}{2}})>0$.

(ii)By choosing a function $\varphi\in D^{1,2}(\mathbb R^{3}), \varphi\geq0, \neq0$ and with $\sup\varphi\subset\{x:f(x)>0\}$, then for $t>0$ small enough, we have
\begin{equation*}
I(t\varphi)\leq\frac{a}{2}t^{2}\|\varphi\|^{2}+\frac{b}{4}t^{4}\|\varphi\|^{4}-\frac{\lambda}{q}t^{q}\int_{\mathbb R^{3}}f(x)|\varphi|^{q}dx-\frac{t^{6}}{6}|\varphi|_{6}^{6}<0.
\end{equation*}
This completes the proof.\quad$\square$\\
Now we are in a position to prove Theorem \ref{thm 1.1}.

\textbf{Proof of Theorem \ref{thm 1.1}} It follows from Lemma \ref{lem 3.1} that
\begin{equation*}
c_{1}=\inf_{u\in\overline{B}_{\alpha}}I(u)<0.
\end{equation*}
By Ekeland's variational principle \cite{MWM}, there exists a minimizing sequence $\{u_{n}\}\subset\overline{B}_{\alpha}$ such that $I(u_{n})\rightarrow c_{1}$ and $I'(u_{n})\rightarrow0$ as $n\rightarrow\infty$. Then, there exists $u_{1}\in D^{1,2}(\mathbb R^{3})$ such that $u_{n}\rightharpoonup u_{1}$ as $n\rightarrow\infty$ in $D^{1,2}(\mathbb R^{3})$. By the fact that $\overline{B}_{\alpha}$ is closed and convex, thus $u_{1}\in\overline{B}_{\alpha}$, then $c_{1}\leq I(u)$. It follows from the H$\ddot{o}$lder inequality that
\begin{equation}\label{c2}
\begin{split}
c_{1}&\leq I(u_{1})\\
&=I(u_{1})-\frac{1}{4}\liminf_{n\rightarrow\infty}(I'(u_{n}), u_{1})\\
&=\frac{a}{4}\|u_{1}\|^{2}+\frac{b}{4}\|u_{1}\|^{2}(\|u_{1}\|^{2}-\liminf_{n\rightarrow\infty}\|u_{n}\|^{2})
+(\frac{1}{4}-\frac{1}{q})\lambda\int_{\mathbb R^{3}}f(x)|u_{1}|^{q}dx+\frac{1}{12}|u_{1}|^{6}_{6}\\
&\leq\liminf_{n\rightarrow\infty}\big{(}\frac{a}{4}\|u_{n}\|^{2}+(\frac{1}{4}-\frac{1}{q})\lambda\int_{\mathbb R^{3}}f(x)|u_{n}|^{q}dx+\frac{1}{12}|u_{n}|^{6}_{6}\big{)}\\
&=\liminf_{n\rightarrow\infty}\big{(}I(u_{n})-\frac{1}{4}(I'(u_{n}), u_{n})\big{)}=c_{1}.
\end{split}
\end{equation}
From above we can deduce that $u_{n}\rightarrow u_{1}$ in $D^{1,2}(\mathbb R^{3})$. Thus the functional $I$ satisfies the $(PS)_{c_{1}}$ condition and the functional $I$ achieves a minimum $u_{1}$ at an interior point of  $B_{\alpha}$. Since $I(u_{1})=I(|u_{1}|)$ we may assume that $u_{1}\geq0$ and by the maximum principle we have $u_{1}>0$ on $\mathbb R^{3}$. This completes the proof.\quad$\square$\\

\section {Proof of Theorem 1.2}
First, we prove the following mountain-pass geometry of functional $I$.
\begin{lemma}\label{lem 4.1}(Mountain pass Geometry) The functional $I$ satisfies the following conditions:
\begin{enumerate}
\item [$(i)$] There exists $\lambda_{1}>0$ such that if $\lambda\in(0,\lambda_{1})$, then there exist $\alpha>0$ and $\rho>0$ such that
\begin{equation*}
I(u)|_{B_{\alpha}}\geq\rho>0.
\end{equation*}
\item [$(ii)$] There exists $e\in D^{1,2}(\mathbb R^{3})$ with $\|e\|>\alpha$ such that $I(e)<0$.
\end{enumerate}
\end{lemma}
\textbf{proof}(i) It directly follows from Lemma \ref{lem 2.1}.

(ii) Note that
\begin{equation*}
I(tU_{\epsilon})=\frac{a}{2}t^{2}\|U_{\epsilon}\|^{2}+\frac{b}{4}t^{4}\|U_{\epsilon}\|^{4}-\frac{t^{q}}{q}\lambda\int_{\mathbb R^{3}}f(x)|U_{\epsilon}|^{q}dx-\frac{t^{6}}{6}|U_{\epsilon}|_{6}^{6}.
\end{equation*}
Then, there exists $t_{0}>0$ sufficient large such that $\|tU_{\epsilon}\|>\rho$, $I(tU_{\epsilon})<0$, where $U_{\epsilon}$ defined in \eqref{b3}. \quad$\square$\\

Therefore, by using the Ambrosetti-Rabinowitz Mountain Pass Theorem without $(PS)_{c_{2}}$ condition(see \cite{MWM}), it follows that there exists a $(PS)_{c_{2}}$ sequence $\{u_{n}\}\subset D^{1,2}(\mathbb R^{3})$ such that
\begin{equation*}
I(u_{n})\rightarrow c_{2}=\inf_{\gamma\in\Gamma}\max_{t\in[0,1]}I_{\lambda}(\gamma(t))\quad \text{and}\quad I'(u_{n})\rightarrow0,
\end{equation*}
where
\begin{equation*}
\Gamma=\{\gamma\in C([0,1], D^{1,2}(\mathbb R^{3})): \gamma(0)=0, I(\gamma(1))<0\}.
\end{equation*}
\begin{lemma}\label{lem 4.2}
The $(PS)_{c_{2}}$ sequence $\{u_{n}\}$ is bounded in $D^{1,2}(\mathbb R^{3})$.
\end{lemma}
\textbf{proof} By the H$\ddot{o}$lder inequality and Lemma \ref{lem 2.1}, we have
\begin{equation*}
\begin{split}
&c_{2}+1+\|u_{n}\|\\
\geq&I(u_{n})-\frac{1}{6}(I'(u_{n}), u_{n})\\
=&\frac{a}{3}\|u_{n}\|^{2}+\frac{b}{12}\|u_{n}\|^{4}-(\frac{1}{q}-\frac{1}{6})\lambda\int_{\mathbb R^{3}}f(x)|u_{n}|^{q}dx\\
\geq&\frac{a}{3}\|u_{n}\|^{2}+\frac{b}{12}\|u_{n}\|^{4}-(\frac{1}{q}-\frac{1}{6})\lambda|f|_{q^{\ast}}|u_{n}|_{6}^{q}\\
\geq&\frac{a}{3}\|u_{n}\|^{2}+\frac{b}{12}\|u_{n}\|^{4}-(\frac{1}{q}-\frac{1}{6})\lambda\frac{1}{S^{\frac{q}{2}}}|f|_{q^{\ast}}\|u_{n}\|^{q}.\\
\end{split}
\end{equation*}
Then $\{u_{n}\}$ is bounded in $D^{1,2}(\mathbb R^{3})$.\quad$\square$\\

\begin{lemma}\label{lem 4.3}
\it If $\{u_{n}\}\subset D^{1,2}(\mathbb R^{3})$ is a bounded $(PS)_{c_{2}}$ sequence of $I$ and $c_{2}<\Lambda-C\lambda^{\frac{2}{2-q}}$, then $\{u_{n}\}$ has a strongly convergent subsequence in $D^{1,2}(\mathbb R^{3})$, where $\Lambda=\frac{1}{4}abS^{3}+\frac{1}{24}b^{3}S^{6}+\frac{1}{24}(b^{2}S^{4}+4aS)^{\frac{3}{2}}$ and $C=\frac{2-q}{2}(\frac{1}{q}-\frac{1}{6})^{\frac{2}{2-q}}|f|_{q^{\ast}}^{\frac{2}{2-q}}(\frac{2q}{aS})^{\frac{q}{2-q}}$.
\end{lemma}
\textbf{proof} By the concentration compactness lemma by P.L. Lions \cite{MWM}, up to a subsequence, there exist an at most countable set $\Gamma$, points $\{a_{k}\}_{k\in\Gamma}\subset\mathbb R^{3}$ and values $\{\eta_{k}\}_{k\in\Gamma}, \{\nu_{k}\}_{k\in\Gamma}\in\mathbb R^{+}$ such that
\begin{equation}\label{d1}
\begin{cases}
|\nabla u_{n}|^{2}\rightharpoonup d\eta\geq|\nabla u|^{2}+\Sigma_{k\in\Gamma}\eta_{k}\delta_{a_{k}},\\
|u_{n}|^{6}\rightharpoonup d\nu=|u|^{2^{\ast}}+\Sigma_{k\in\Gamma}\nu_{k}\delta_{a_{k}},\\
\end{cases}
\end{equation}
where $\delta_{a_{k}}$ is the Dirac delta measure concentrated $a_{k}$. Moreover,
\begin{equation}\label{d2}
\nu_{k}\leq\eta_{k}^{3}S^{-3}.
\end{equation}
In the following, we prove that $\Gamma=\varnothing$. Arguing by contradiction, fix $k\in\Gamma$, for $\epsilon>0$, assume that $\psi_{\epsilon}^{k}\in C_{0}^{\infty}(\mathbb R^{3}, [0, 1])$ such that
\begin{equation*}
\begin{cases}
\psi_{\epsilon}^{k}=1,\ &for~~|x-a_{k}|\leq\frac{\epsilon}{2},\\
\psi_{\epsilon}^{k}=0,\ &for~~|x-a_{k}|\geq\epsilon,\\
|\nabla \psi_{\epsilon}^{k}|\leq\frac{3}{\epsilon},\ &in~~\mathbb R^{3}.\\
\end{cases}
\end{equation*}
Since $\{\psi_{\epsilon}^{k}u_{n}\}$ is bounded in $D^{1,2}(\mathbb R^{3})$, we have
\begin{equation*}
(I'(u_{n}), \psi_{\epsilon}^{k}u_{n})\rightarrow0,
\end{equation*}
i.e.
\begin{equation}\label{d3}
(a+b\|u_{n}\|^{2})(\int_{\mathbb R^{3}}u_{n}\triangledown u_{n}\triangledown \psi_{\epsilon}^{k}dx+\int_{\mathbb R^{3}}|\triangledown u_{n}|^{2} \psi_{\epsilon}^{k}dx)=\lambda\int_{\mathbb R^{3}}f(x)|u_{n}|^{q}\psi_{\epsilon}^{k}dx+\int_{\mathbb R^{3}}|u_{n}|^{6}\psi_{\epsilon}^{k}dx+o(1).
\end{equation}
 It follows from the boundedness of $\{u_{n}\}$ in $D^{1,2}(\mathbb R^{3})$ and the H$\ddot{o}$lder inequality that
\begin{equation}\label{d4}
\begin{split}
&\lim_{\epsilon\rightarrow0}\limsup_{n\rightarrow\infty}(a+b\|u_{n}\|^{2})\int_{\mathbb R^{3}}u_{n}\triangledown u_{n}\triangledown \psi_{\epsilon}^{k}dx\\
\leq&\lim_{\epsilon\rightarrow0}\limsup_{n\rightarrow\infty}C_{1}(\int_{B_{\epsilon}(a_{k})}|\triangledown u_{n}|^{2}dx)^{\frac{1}{2}}(\int_{B_{\epsilon}(a_{k})}|\triangledown \psi_{\epsilon}^{k}|^{2}|u_{n}|^{2}dx)^{\frac{1}{2}}\\
\leq&\lim_{\epsilon\rightarrow0}C_{2}(\int_{B_{\epsilon}(a_{k})}|\triangledown \psi_{\epsilon}^{k}|^{2}|u|^{2}dx)^{\frac{1}{2}}\\
\leq&\lim_{\epsilon\rightarrow0}C_{3}(\int_{B_{\epsilon}(a_{k})}|\triangledown \psi_{\epsilon}^{k}|^{3}dx)^{\frac{1}{3}}(\int_{B_{\epsilon}(a_{k})}|u|^{6}dx)^{\frac{1}{6}}\\
\leq&\lim_{\epsilon\rightarrow0}C_{4}(\int_{B_{\epsilon}(a_{k})}|u|^{6}dx)^{\frac{1}{6}}\\
=&0,
\end{split}
\end{equation}
and
\begin{equation}\label{d5}
\begin{split}
&\lim_{\epsilon\rightarrow0}\limsup_{n\rightarrow\infty}\int_{\mathbb R^{3}}f(x)|u_{n}|^{q}\psi_{\epsilon}^{k}dx\\
=&\lim_{\epsilon\rightarrow0}\int_{B_{\epsilon}(a_{k})}f(x)|u|^{q}\psi_{\epsilon}^{k}dx\\
=&0.
\end{split}
\end{equation}
From \eqref{b5}, we have
\begin{equation}\label{d6}
\begin{split}
&\lim_{\epsilon\rightarrow0}\limsup_{n\rightarrow\infty}(a+b\|u_{n}\|^{2})\int_{\mathbb R^{3}}|\triangledown u_{n}|^{2} \psi_{\epsilon}^{k}dx\\
\geq&\lim_{\epsilon\rightarrow0}\limsup_{n\rightarrow\infty}a\int_{\mathbb R^{3}}|\triangledown u_{n}|^{2} \psi_{\epsilon}^{k}dx+\limsup_{n\rightarrow\infty}b(\int_{\mathbb R^{3}}|\triangledown u_{n}|^{2} \psi_{\epsilon}^{k}dx)^{2}\\
\geq&a\eta_{k}+b\eta_{k}^{2},
\end{split}
\end{equation}
and
\begin{equation}\label{d7}
\begin{split}
&\lim_{\epsilon\rightarrow0}\limsup_{n\rightarrow\infty}\int_{\mathbb R^{3}}|u_{n}|^{6}\psi_{\epsilon}^{k}dx\\
=&\lim_{\epsilon\rightarrow0}\int_{\mathbb R^{3}}|u|^{6}\psi_{\epsilon}^{k}dx+\nu_{k}\\
=&\nu_{k}.
\end{split}
\end{equation}
From \eqref{d3}-\eqref{d7}, we have
\begin{equation}\label{d8}
\nu_{k}\geq a\eta_{k}+b\eta_{k}^{2}.
\end{equation}
Combining with \eqref{d2}, we can deduce that
\begin{equation}\label{d9}
\eta_{k}\geq\frac{b+\sqrt{b^{2}+4aS^{-3}}}{2S^{-3}}.
\end{equation}

For $R>0$, assume that $\phi_{R}\in C_{0}^{\infty}(\mathbb R^{3}, [0, 1])$ such that
\begin{equation*}
\begin{cases}
\phi_{R}=1,\ &for~~|x|<R,\\
\phi_{R}=0,\ &for~~|x|\geq 2R,\\
|\nabla \phi_{R}|\leq\frac{2}{R},\ &in~~\mathbb R^{3}.\\
\end{cases}
\end{equation*}
Then, by Lemma \ref{lem 2.2}, we have
\begin{equation}\label{d10}
\begin{split}
c_{2}=&\lim_{n\rightarrow\infty}(I(u_{n})-\frac{1}{4}(I'(u_{n}), u_{n})\\
=&\lim_{n\rightarrow\infty}(\frac{a}{4}\|u_{n}\|^{2}+\frac{1}{12}\int_{\mathbb R^{3}}|u_{n}|^{6}dx-(\frac{1}{q}-\frac{1}{4})\lambda\int_{\mathbb R^{3}}f(x)|u_{n}|^{q}dx\\
\geq&\lim_{R\rightarrow\infty}\lim_{n\rightarrow\infty}(\frac{a}{4}\int_{\mathbb R^{3}}|\nabla u_{n}|^{2}\phi_{R}dx+\frac{1}{12}\int_{\mathbb R^{3}}|u_{n}|^{6}\phi_{R}dx-(\frac{1}{q}-\frac{1}{4})\lambda\int_{\mathbb R^{3}}f(x)|u_{n}|^{q}dx\\
\geq&\lim_{R\rightarrow\infty}(\frac{a}{4}\int_{\mathbb R^{3}}|\nabla u_{n}|^{2}\phi_{R}dx+\frac{a}{4}\eta_{k}+\frac{1}{12}\int_{\mathbb R^{3}}|u_{n}|^{6}\phi_{R}dx+\frac{1}{12}\nu_{k}-(\frac{1}{q}-\frac{1}{4})\lambda\int_{\mathbb R^{3}}f(x)|u|^{q}dx\\
\geq&\frac{a}{4}\int_{\mathbb R^{3}}|\nabla u|^{2}dx+\frac{a}{4}\eta_{k}+\frac{1}{12}\int_{\mathbb R^{3}}|u_{n}|^{6}\phi_{R}d+\frac{1}{12}\nu_{k}-(\frac{1}{q}-\frac{1}{4})\lambda\int_{\mathbb R^{3}}f(x)|u|^{q}dx\\
\geq&\frac{aS}{4}|u|_{6}^{2}+\frac{a}{4}\eta_{k}+\frac{1}{12}\nu_{k}-(\frac{1}{q}-\frac{1}{4})\lambda\int_{\mathbb R^{3}}f(x)|u|^{q}dx.\\
\end{split}
\end{equation}
By \eqref{d8} and \eqref{d9}, we have
\begin{equation}\label{d11}
\frac{a}{4}\eta_{k}+\frac{1}{12}\nu_{k}\geq\Lambda.
\end{equation}
In order to estimate $\frac{aS}{4}|u|_{6}^{2}-(\frac{1}{q}-\frac{1}{4})\int_{\mathbb R^{3}}f(x)|u|^{q}dx$, we observe that the function
\begin{equation}\label{d12}
\begin{split}
t\mapsto&\frac{aS}{4}t^{2}-(\frac{1}{q}-\frac{1}{4})\lambda\int_{\mathbb R^{3}}f(x)|u|^{q}dx\\
\geq&\frac{aS}{4}t^{2}-(\frac{1}{q}-\frac{1}{4})\lambda|f|_{q^{\ast}}t^{q}:=\widetilde{f}(t)\\
\end{split}
\end{equation}
achieves its minimum on $(0, \infty)$ at a point  $t_{1}$, $\min_{t\geq0}\widetilde{f}(t)=\widetilde{f}(t_{1})=-C\lambda^{\frac{2}{2-q}}$, where \begin{equation*}
t_{1}=\frac{2q}{aS}(\frac{1}{q}-\frac{1}{4})|f|_{q^{\ast}}^{\frac{1}{2-q}} \quad\text{and}\quad C=\frac{2-q}{2}(\frac{1}{q}-\frac{1}{6})^{\frac{2}{2-q}}|f|_{q^{\ast}}^{\frac{2}{2-q}}(\frac{2q}{aS})^{\frac{q}{2-q}}.
\end{equation*}
From \eqref{d10}, \eqref{d11} and \eqref{d12}, we can deduce that $c_{2}\geq\Lambda-C\lambda^{\frac{2}{2-q}}$, which is a contradiction. Thus $\Gamma=\varnothing$.

For $R>0$, define
\begin{equation}\label{d13}
\eta_{\infty}=\lim_{R\rightarrow\infty}\limsup_{n\rightarrow\infty}\int_{|x|>R}|\triangledown u_{n}|^{2}dx
\end{equation}
and
\begin{equation}\label{d14}
\nu_{\infty}=\lim_{R\rightarrow\infty}\limsup_{n\rightarrow\infty}\int_{|x|>R}|u_{n}|^{6}dx.
\end{equation}
Then,
\begin{equation}\label{d15}
\limsup_{n\rightarrow\infty}\int_{\mathbb R^{3}}|\triangledown u_{n}|^{2}dx=\int_{\mathbb R^{3}}d\eta+\eta_{\infty}\quad\text{and}\quad
\limsup_{n\rightarrow\infty}\int_{\mathbb R^{3}}|u_{n}|^{6}dx=\int_{\mathbb R^{3}}d\nu+\nu_{\infty}.
\end{equation}
Moreover,
\begin{equation}\label{d16}
\nu_{\infty}\leq\eta_{\infty}^{3}S^{-3}.
\end{equation}
Assume that $\chi_{R}\in C_{0}^{\infty}(\mathbb R^{3}, [0, 1])$ such that
\begin{equation*}
\begin{cases}
\chi_{R}(x)=0,\ &for~~|x|<\frac{R}{2},\\
\chi_{R}(x)=1,\ &for~~|x|>R,\\
|\triangledown \chi_{R}(x)|<\frac{3}{R},\ &in~~\mathbb R^{3}.\\
\end{cases}
\end{equation*}
Since $\{\chi_{R}u_{n}\}$ is bounded in $D^{1,2}(\mathbb R^{3})$, we have
\begin{equation*}
(I'(u_{n}), \chi_{R}u_{n})\rightarrow0,
\end{equation*}
i.e.
\begin{equation}\label{d17}
(a+b\|u_{n}\|^{2})(\int_{\mathbb R^{3}}u_{n}\triangledown u_{n}\triangledown \chi_{R}dx+\int_{\mathbb R^{3}}|\triangledown u_{n}|^{2} \chi_{R}dx)=\lambda\int_{\mathbb R^{3}}f(x)|u_{n}|^{q}\chi_{R}dx+\int_{\mathbb R^{3}}|u_{n}|^{6}\chi_{R}dx+o(1).
\end{equation}
 It follows from the boundedness of $\{u_{n}\}$ in $D^{1,2}(\mathbb R^{3})$ and the H$\ddot{o}$lder inequality that
\begin{equation}\label{d18}
\begin{split}
&\lim_{R\rightarrow\infty}\limsup_{n\rightarrow\infty}(a+b\|u_{n}\|^{2})\int_{\mathbb R^{3}}u_{n}\triangledown u_{n}\triangledown \chi_{R}dx\\
\leq&\lim_{R\rightarrow\infty}\limsup_{n\rightarrow\infty}(\int_{\frac{R}{2}\leq|x|\leq R}|\triangledown u_{n}|^{2}dx)^{\frac{1}{2}}(\int_{\frac{R}{2}\leq|x|\leq R}|\triangledown\chi_{R}|^{2}|u_{n}|^{2}dx)^{\frac{1}{2}}\\
\leq&\lim_{R\rightarrow\infty}C_{2}(\int_{\frac{R}{2}\leq|x|\leq R}|\triangledown \chi_{R}^{k}|^{2}|u|^{2}dx)^{\frac{1}{2}}\\
\leq&\lim_{R\rightarrow\infty}C_{3}(\int_{\frac{R}{2}\leq|x|\leq R}|\triangledown \chi_{R}^{k}|^{3}dx)^{\frac{1}{3}}(\int_{\frac{R}{2}\leq|x|\leq R}|u|^{6}dx)^{\frac{1}{6}}\\
\leq&\lim_{R\rightarrow\infty}C_{4}(\int_{\frac{R}{2}\leq|x|\leq R}|u|^{6}dx)^{\frac{1}{6}}\\
=&0,
\end{split}
\end{equation}
and
\begin{equation}\label{d19}
\begin{split}
&\lim_{R\rightarrow\infty}\limsup_{n\rightarrow\infty}\int_{\mathbb R^{3}}f(x)|u_{n}|^{q}\chi_{R}dx\\
=&\lim_{R\rightarrow\infty}\int_{\mathbb R^{3}}f(x)|u|^{q}\chi_{R}dx\\
=&0.
\end{split}
\end{equation}
From \eqref{d13} and \eqref{d14}, we have
\begin{equation}\label{d20}
\begin{split}
&\lim_{R\rightarrow\infty}\limsup_{n\rightarrow\infty}(a+b\|u_{n}\|^{2})\int_{\mathbb R^{3}}|\triangledown u_{n}|^{2}\chi_{R}dx\\
\geq&\lim_{R\rightarrow\infty}\limsup_{n\rightarrow\infty}a\int_{\mathbb R^{3}}|\triangledown u_{n}|^{2} \chi_{R}dx+\lim_{R\rightarrow\infty}\limsup_{n\rightarrow\infty}b(\int_{\mathbb R^{3}}|\triangledown u_{n}|^{2}\chi_{R}dx)^{2}\\
\geq&\lim_{R\rightarrow\infty}\limsup_{n\rightarrow\infty}a\int_{|x|>R}|\triangledown u_{n}|^{2}dx+\lim_{R\rightarrow\infty}\limsup_{n\rightarrow\infty}b(\int_{|x|>R}|\triangledown u_{n}|^{2}dx)^{2}\\
=&a\eta_{\infty}+b\eta_{\infty}^{2}
\end{split}
\end{equation}
and
\begin{equation}\label{d21}
\begin{split}
&\lim_{R\rightarrow\infty}\limsup_{n\rightarrow\infty}\int_{\mathbb R^{3}}|u_{n}|^{6}\chi_{R}dx\\
=&\lim_{R\rightarrow\infty}\limsup_{n\rightarrow\infty}\int_{|x|>\frac{R}{2}}|u_{n}|^{6}\chi_{R}dx\\
\leq&\lim_{R\rightarrow\infty}\limsup_{n\rightarrow\infty}\int_{|x|>\frac{R}{2}}|u_{n}|^{6}dx\\
=&\nu_{\infty}.
\end{split}
\end{equation}
From \eqref{d17}-\eqref{d21}, we have
\begin{equation}\label{d22}
\nu_{\infty}\geq a\eta_{\infty}+b\eta_{\infty}^{2}.
\end{equation}
From $\eta_{\infty}, \nu_{\infty}\geq0$, \eqref{d16} and \eqref{d22}, we can see that $\eta_{\infty}=0$ if and only if $\nu_{\infty}=0$. In the following, we assume that $\eta_{\infty}\neq0$. From \eqref{d16} and \eqref{d22}, we can deduce that
\begin{equation}\label{d23}
\eta_{\infty}\geq\frac{b+\sqrt{b^{2}+4aS^{-3}}}{2S^{-3}}.
\end{equation}
Then, it follows from \eqref{d12}-\eqref{d14}, \eqref{d22} and \eqref{d23}, we have
\begin{equation}\label{d24}
\begin{split}
c_{2}=&\lim_{n\rightarrow\infty}(I(u_{n})-\frac{1}{4}(I'(u_{n}), u_{n})\\
=&\lim_{n\rightarrow\infty}(\frac{a}{4}\|u_{n}\|^{2}+\frac{1}{12}\int_{\mathbb R^{3}}|u_{n}|^{6}dx-(\frac{1}{q}-\frac{1}{4})\lambda\int_{\mathbb R^{3}}f(x)|u_{n}|^{q}dx\\
\geq&\frac{a}{4}\int_{\mathbb R^{3}}d\eta+\frac{a}{4}\eta_{\infty}+\frac{1}{12}\int_{\mathbb R^{3}}d\nu+\frac{1}{12}\nu_{\infty}-(\frac{1}{q}-\frac{1}{4})\lambda\int_{\mathbb R^{3}}f(x)|u|^{q}dx\\
\geq&\frac{a}{4}\int_{\mathbb R^{3}}|\nabla u|^{2}dx+\frac{a}{4}\eta_{\infty}+\frac{1}{12}\nu_{\infty}-(\frac{1}{q}-\frac{1}{4})\lambda\int_{\mathbb R^{3}}f(x)|u|^{q}dx\\
\geq&\frac{aS}{4}|u|_{6}^{2}+\frac{a}{4}\eta_{\infty}+\frac{1}{12}\nu_{\infty}-(\frac{1}{q}-\frac{1}{4})\lambda\int_{\mathbb R^{3}}f(x)|u|^{q}dx\\
\geq&\Lambda-C\lambda^{\frac{2}{2-q}},
\end{split}
\end{equation}
which contradict with the assumption of $c_{2}\leq\Lambda-C\lambda^{\frac{2}{2-q}}$. Thus $\eta_{\infty}=\nu_{\infty}=0$.
From $\Gamma=\varnothing$ and \eqref{d14}, we can obtain that
\begin{equation*}
\limsup_{n\rightarrow\infty}\int_{\mathbb R^{3}}|u_{n}|^{6}dx=\int_{\mathbb R^{3}}|u|^{6}dx.
\end{equation*}
Combining with Fatou Lemma, we have
\begin{equation*}
\int_{\mathbb R^{3}}|u|^{6}dx\leq\liminf_{n\rightarrow\infty}\int_{\mathbb R^{3}}|u_{n}|^{6}dx\leq\limsup_{n\rightarrow\infty}\int_{\mathbb R^{3}}|u_{n}|^{6}dx.
\end{equation*}
Thus, $\lim_{n\rightarrow\infty}\int_{\mathbb R^{3}}|u_{n}|^{6}dx=\int_{\mathbb R^{3}}|u|^{6}dx$.

In the following that we prove that $u_{n}\rightarrow u$ in $D^{1,2}(\mathbb R^{3})$. Assume that $\lim_{n\rightarrow\infty}\|u_{n}\|^{2}=d^{2}$, it sufficient to show that $d^{2}=\int_{\mathbb R^{3}}|\triangledown u|^{2}dx$. Indeed,
\begin{equation*}
\begin{split}
0=&\lim_{n\rightarrow\infty}(I'(u_{n}), u_{n}-u)=\lim_{n\rightarrow\infty}I'(u_{n})u_{n}-\lim_{n\rightarrow\infty}I'(u_{n})u\\
=&\lim_{n\rightarrow\infty}(a+b\|u_{n}\|^{2})\|u_{n}\|^{2}-\int_{\mathbb R^{3}}|u_{n}|^{6}dx-\lambda\int_{\mathbb R^{3}}f(x)|u_{n}|^{q}dx\\
-&\lim_{n\rightarrow\infty}[(a+b\|u_{n}\|^{2})\int_{\mathbb R^{3}}\nabla u_{n}\nabla udx-\int_{\mathbb R^{3}}|u_{n}|^{4}u_{n}udx-\lambda\int_{\mathbb R^{3}}f(x)|u_{n}|^{q-2}u_{n}udx]\\
=&(a+bd^{2})(d^{2}-\int_{\mathbb R^{3}}|\nabla u|^{2}dx).
\end{split}
\end{equation*}
Thus, $u_{n}\rightarrow u$ in $D^{1,2}(\mathbb R^{3})$. This shows that the functional $I$ satisfies the $(PS)_{c_{2}}$ 
condition for $c_{2}<\Lambda-C\lambda^{\frac{2}{2-q}}$.\quad$\square$\\

In the following, we estimate the energy $c_{2}$.
\begin{lemma}\label{lem 4.4}
There exists $\lambda_{2}\in(0, \lambda_{1}]$ such that for any $\lambda\in(0, \lambda_{2})$, we have
\begin{equation*}
c_{2}\leq\sup_{t\geq0}I(tU_{\epsilon})<\Lambda-C\lambda^{\frac{2}{2-q}}.
\end{equation*}
\end{lemma}
\textbf{proof} It follows from \eqref{b4} that
\begin{equation*}
\begin{split}
I(tU_{\epsilon})=&\frac{at^{2}}{2}\|U_{\epsilon}\|^{2}+\frac{bt^{4}}{4}\|U_{\epsilon}\|^{4}-\frac{t^{6}}{6}|U_{\epsilon}|^{6}_{6}-\frac{t^{q}}{q}\lambda\int_{\mathbb R^{3}}f(x)|U_{\epsilon}|^{q}dx\\
=&\frac{at^{2}}{2}S^{\frac{3}{2}}+\frac{bt^{4}}{4}S^{3}-\frac{t^{6}}{6}S^{\frac{3}{2}}-\frac{t^{q}}{q}\lambda\int_{\mathbb R^{3}}f(x)|U_{\epsilon}|^{q}dx.\\
\end{split}
\end{equation*}
We observe that the function
\begin{equation*}
\begin{split}
t\mapsto&\frac{at^{2}}{2}\|U_{\epsilon}\|+\frac{bt^{4}}{4}\|U_{\epsilon}\|-\frac{t^{6}}{6}|U_{\epsilon}|^{6}_{6}\\
&=\frac{at^{2}}{2}S^{\frac{3}{2}}+\frac{at^{4}}{4}S^{3}-\frac{t^{6}}{6}S^{\frac{3}{2}}:=\widetilde{g}(t)
\end{split}
\end{equation*}
achieves its maximum on $[0, \infty)$ at a point  $t_{2}$, that is, $\max_{t\geq0}\widetilde{g}(t)=\widetilde{g}(t_{2})=\Lambda$,
where $t_{2}^{2}=\frac{bS^{3}+\sqrt{b^{2}S^{6}+4aS^{3}}}{2S\frac{3}{2}}$.
First, we choose $t_{3}\in(0, t_{2})$ and $\lambda_{3}>0$ small enough such that
\begin{equation*}
I(tU_{\epsilon})<\Lambda-C\lambda^{\frac{2}{2-q}}
\end{equation*}
for $0\leq t\leq t_{3}$ and $0<\lambda\leq\lambda_{3}$. Here $\lambda_{3}$ is chosen so that
\begin{equation*}
\Lambda-C\lambda^{\frac{2}{2-q}}>0
\end{equation*}
for all $0<\lambda\leq\lambda_{3}$. To estimate $I(tU_{\epsilon})$ for $t\geq t_{3}$,
by $1\leq q<2$, we can choose $\lambda_{2}\in(0, \lambda_{3}]$ such that for any $\lambda\in(0, \lambda_{2}]$, we have
\begin{equation*}
\frac{t_{3}^{q}}{q}\lambda\int_{\mathbb R^{3}}f(x)|U_{\epsilon}|^{q}dx>C\lambda^{\frac{2}{2-q}}.
\end{equation*}
Then, for any $\lambda\in(0, \lambda_{2})$, we have
\begin{equation*}
\begin{split}
\sup_{t\geq t_{3}}I(tU_{\epsilon})&\leq\sup_{t\geq t_{3}}(\widetilde{g}(t)-\frac{t_{3}^{q}}{q}\lambda\int_{\mathbb R^{3}}f(x)|U_{\epsilon}|^{q}dx)\\
&<\widetilde{g}(t_{1})-C\lambda^{\frac{2}{2-q}}\\
&=\Lambda-C\lambda^{\frac{2}{2-q}}.
\end{split}
\end{equation*}
Thus, we complete the proof.\quad$\square$\\

\textbf{Proof of Theorem \ref{thm 1.2}} The argument of Theorem \ref{thm 1.1} shows that there exist $\lambda_{1}>0$ such that for 
each $\lambda\in(0, \lambda_{1})$ problem \eqref{a1} has a solution which is a local minimum of $I$.
From Lemmas \ref{lem 4.2}, \ref{lem 4.3} and \ref{lem 4.4}, there exist $\lambda_{2}\leq\lambda_{1}$ such that for each $\lambda\in(0, \lambda_{2})$,
problem \eqref{a1} has a solution which is a mountain pass solution.\quad$\square$\\


\begin{thebibliography}{99}

\bibitem{AAHBGC} A. Ambrosetti, H. Brezis, G. Cerami, Combined effects of concave and convex nonlinearities in some elliptic problems, J. Funct. Anal. 122 (1994) 519-543.
\bibitem{AAAZZ}A. Azzollini, P. d$'$Avenia, A. Pomponio, Multiple critical points for a class of nonlinear functionals, Ann. Mat. Pura Appl. 190 (2011) 507-523.

\bibitem{TBTBT} T. Bartsch, Infinitely many solutions of a symmetric Dirichlet problem, Nonlinear Anal. TMA, 20 (1993) 1205-1206.

\bibitem{HBPLL} H. Berestycki, P.L. Lions, Nonlinear scalar field equations, I existence of a ground state, Archive for Rational Mechanics and Analysis, 82(4) (1983) 313-345.

 \bibitem{BARTS}T. Bartsch, Z.Q. Wang, Existence and multiple results for some superlinear elliptic problems on $\mathbb{R}^N$, Commun. Partial Differ. Equ. 20 (1995) 1725-1741.

\bibitem{KBTW} K.J. Brown, T.F. Wu, A fibering map approach to a semilinear elliptic booundary value problem, Electron. J. Differential Equations 69 (2007) 1-9.

\bibitem{KBYZ} K.J. Brown, Y.P. Zhang, The Nehari manifold for a semilinear elliptic equation with a sign-changing weight function, J. Differential Equations 193 (2003) 481-499.

\bibitem{Chabrowski} J. Chabrowski, P. Drabek, On positive solutions of nonlinear elliptic equations involving concave and critical nonlinearities, Studia mathematica 151 (2001) 67-85.

\bibitem{KCC} K.C. Chang, Methods in Nonlinear Analysis, Springer, 2005.

\bibitem{CCYKTW} C.Y. Chen, Y.C. Kuo, T.F. Wu, The Nehari manifold for a Kirchhoff type problem involving sign-changing weight functions, J. Differential Equations 250 (2011) 1876-1908.

\bibitem{CLLM} S.J. Chen, L. Li, Multiple solutions for the nonhomogeneous Kirchhoff equation on $R^{N}$, Nonlinear Anal. Real World Appl. 14 (2013) 1477-1486.

\bibitem{BCBC} B.T. Cheng, New existence and multiplicity of nontrivial solution for nonlocal elliptic Kirchhoff type problems, J. Math. Anal. Appl. 394 (2012) 488-495.

\bibitem{BCXW} B.T. Cheng, X. Wu, Existence results of positive solutions of Kirchhoff type problems, Nonlinear Anal. 71 (2009) 4883-4892.

\bibitem{DGW} G.W. Dai, R.F. Hao, Existence of solutions for a {$p(x)$}-{K}irchhoff-type equation, J. Math. Anal. Appl. 1 (2009) 275--284.

\bibitem{YDYD} Y.H. Ding, Variational Methods for Strongly Indefinite Problems, World Scientific Press. 2008.

\bibitem{PDSP} P. Dr$\acute{a}$bek, S. Poho$\breve{z}$ev, Positive solutions for the p-Laplacian: application of the fibering method, Proc. Roy. Soc. Edinburgh Sect. A 127 (1997) 703-726.

\bibitem{Fan} H.N. Fan, Multiple positive solutions for a class of Kirchhoff type problems involving critical Soblev exponents, J. Math. Anal. Appl. 431 (2015) 150-168.

\bibitem{GMFG} Giovany M. Figueiredo,  Existence of a positive solution for a Kirchhoff problem type with critical growth via truncation argument, J. Math. Anal. Appl. 401 (2013) 706-713.

\bibitem{XHWZ} X.M. He, W.M. Zou, Multiplicity of solutions for a class of Kirchhoff type problems, Acta Math. Appl. Sin. Engl. Ser. 26 (2010) 387-394.

\bibitem{XHWWZ} X.M. He, W.M. Zou, Existence and concentration behavior of positive solutions for a Kirchhoff equation in $R^{3}$, J. Differential Equations 2 (2012) 1813-1834.

\bibitem{JJXW}J.H. Jin, X. Wu, Infinitely many radial solution for Kirchhoff-type problems in $R^{N}$, J. Math. Anal. Appl. 369 (2010) 564-574.

\bibitem{GGKKM} G. Kirchhoff, Mechanik, Teubner, Leipzig, 1883.

\bibitem{JLLLL} J.L. Lions, On some questions in boundary value problems of mathematical physics, in: Contemporary Developments in Continuum Mechanics and Partial Differential Equations, Proceedings of International Symposium, Inst. Mat., Univ. Fed. Riode Janeiro, Rio de Janeiro, 1977, in: North-Holland Math. Stud., vol. 30, North-Holland, Amsterdam, 1978, pp. 284-346.

\bibitem{GGBBL} G.B. Li, H.Y. Ye, Existence of positive ground state solutions for the nonlinear Kirchhoff type equations in $\mathbb{R}^3$, J. Differential Equations, 257 (2014) 566-600.

\bibitem{LFYL} Y.H. Li, F.Y. Li, J.P. Shi,  Existence of a positive solution to Kirchhoff type problems   without compactness conditions, J. Differential Equations 253 (2012) 2285-2294.

\bibitem{YLFL} Y.H. Li, F.Y. Li, J.P. Shi, Existence of positive solutions to Kirchhoff type problems with zero mass, J. Math. Anal. Appl. 410 (2014) 361-374.

\bibitem{ZLFL}Z.P. Liang, F.Y. Li, J.P. Shi, Positive solutions to Kirchhoff type equations with nonlinearity having prescribed asymptotic behavior,  Ann. I. H. Poincar$\acute{e}$- AN 31 (2014) 155-167.

\bibitem{Liu}J. Liu, J.F. Liao, C.L. Tang, Positive solutions for Kirchhoff-type equations with critical exponent in $\mathbb{R}^N$, J. Math. Anal. Appl. 429 (2015) 1153-1172.

\bibitem{WLXH} W. Liu, X.M. He, Multiplicity of high energy solutions for superlinear Kirchhoff equations, J.Appl.Math.Comput. 39 (2012) 473-487.

\bibitem{ZLLZQW}Z.L. Liu, Z.Q. Wang, Multiple bound states of nonlinear Schr$\ddot{o}$dinger systems, Comm. Math. Phys. 282 (2008) 721-731.

\bibitem{ZSLSJG} Z.S. Liu, S.J. Guo, On ground states for the Kirchhoff-type problem with a general critical nonlinearity,  J. Math. Anal. Appl. 426 (2015) 267-287.

\bibitem{AMZZ} A.M. Mao, Z.T. Zhang, Sign-changing and multiple solutions of Kirchhoff type problems without the P.S. condition, Nonlinear Anal. 70 (2009) 1275-1287.

\bibitem{JSSL}J. Sun, S.B. Liu, Nontrival solution of Kirchhoff type problems, Appl. Math. Lett. 25 (2012) 500-504.

\bibitem{YZXL} Y.J. Sun, X. Liu, Existence of Positive Solutions for Kirchhoff Type Problems with Critical Exponent, J. Part. Diff. Eq. 25 (2012) 85-96.

\bibitem{HSCT}J.J. Sun, C.L. Tang, Existence and multiplicity of solutions for Kirchhoff type equations, Nonlinear Anal. 74 (2011) 1212-1222.

\bibitem{ASTW} A. Szulkin, T. Weth, The method of Nahari manifold, Boston, 2010, 597-632.



\bibitem{MWM} M. Willem, Minimax Theorems, Progr. Nonlinear Differential Equations Appl., vol. 24, Birkhauser, Basel, 1996.

\bibitem{TWU} T.F. Wu, Multiple positive solutions for a class of concave-convex ellipic problems in $\mathbb{R}^{N}$ involving sign-changing weight, J. Funct. Anal. 258 (2010) 99-131.

\bibitem{XXWW}X. Wu, Existence of nontrivial solutions and high energy solutions for Schroinger-Kirchhoff type equations in $R^{N}$, Nonlinear Anal. Real World Appl. 12 (2011) 1278-1287.

\bibitem{JWLX} J. Wang, L. X. Tian, J.X. Xu, F.B. Zhang, Multiplicity and concentration of positive solutions for a Kirchhoff type problem with critical growth, J. Differential Equations 253 (2012) 2314-2351.

\bibitem{Xie} Q.L. Xie, S.W. Ma, X. Zhang, Bound state solutions of Kirchhoff type problems with critical exponent, J. Differential Equations (2016), http://dx.doi.org/10.1016/j.jde.2016.03.028.

\bibitem{ZTZKP} Z.T. Zhang, K. Perera, Sign changing solutions of Kirchhoff type problems via invariant sets of descent flow, J. Math. Anal. Appl. 317 (2)
(2006) 456-463.

\bibitem{WZMSC} W.M. Zou, M. Schechter, Critical Point Theory and Its Applications, Springer, New York 2006.

\bibitem{HZFZ} H. Zhang, F.B. Zhang, Ground states for the nonlinear Kirchhoff type problems, J. Math. Anal. Appl. 423 (2015) 1671-1692.








\end{thebibliography}
\end{document}